\documentclass[12pt]{article}
\usepackage{amsmath}
\usepackage{amssymb}
\usepackage{vatola}

\def\q{\quad}
\def\qq{\qquad}
\def\mod{\pmod}
\def\t{\text}
\def\f{\frac}
\def\e{\equiv}

\def\a{\alpha}
\def\qtq#1{\q\t{#1}\q}
\def\sls#1#2{(\f{#1}{#2})}
 \def\ls#1#2{\big(\f{#1}{#2}\big)}
\def\Ls#1#2{\Big(\f{#1}{#2}\Big)}
\let \pro=\proclaim
\let \endpro=\endproclaim

\begin{document}
\leftline{The paper will appear in International Journal of Number
Theory.}
\par\q\par\q
\centerline {\bf On the number of representations of $n$ as a}
\centerline {\bf \qq \qq linear combination of four triangular
numbers}

$$\q$$
\centerline{Min Wang$^1$ and Zhi-Hong Sun$^2$}
\par\q\newline
\centerline{$\ ^1$School of Mathematical Sciences, Soochow
University,}
 \centerline{Suzhou, Jiangsu 215006,
 P.R. China}
\centerline{Email: 20144207002@stu.suda.edu.cn}
\par\q\newline
\centerline{$\ ^2$School of Mathematical Sciences, Huaiyin Normal
University,} \centerline{Huaian, Jiangsu 223001, P.R. China}
\centerline{Email: zhihongsun@yahoo.com} \centerline{Homepage:
http://www.hytc.edu.cn/xsjl/szh}

 \abstract{Let $\Bbb Z$ and $\Bbb N$ be the set of integers
 and the set of positive integers, respectively. For
 $a,b,c,d,n\in\Bbb N$ let $t(a,b,c,d;n)$ be the number of
 representations of $n$ by $ax(x-1)/2+by(y-1)/2+cz(z-1)/2
 +dw(w-1)/2$ $(x,y,z,w\in\Bbb Z$). In this paper we obtain explicit
 formulas for $t(a,b,c,d;n)$ in the cases
 $(a,b,c,d)=(1,2,2,4),\ (1,2,4,4),\ (1,1,4,4),\ (1,4,4,4)$,
 $(1,3,9,9),\ (1,1,3,9)$, $(1,3,3,9)$,
 $(1,1,9,9),\ (1,9,9,9)$ and $(1,1,1,9).$
 \par\q
 }
 \endabstract
\let\thefootnote\relax \footnotetext {The second author is the corresponding
author.}

\section*{1. Introduction}
\par\q  Let $\Bbb Z$ and $\Bbb N$ be the set of integers
 and the set of positive integers, respectively.
For  $n \in \Bbb N$ let
$$\sigma(n)=\sum_{d \mid n,d\in\Bbb N}d.$$ For convenience
 we define $\sigma(n)=0$ for $n\notin \Bbb N$. For $a,b,c,d\in\Bbb N$ and $n\in\Bbb N \cup \{0\}$ set
$$N(a,b,c,d;n)=\big|\{(x,y,z,w)\in \Bbb Z^4\ |\ n=ax^2+by^2+cz^2+dw^2
\}\big|$$ and $$t(a,b,c,d;n)=\Big|\Big\{(x,y,z,w)\in \Bbb Z^4\ |\ n\
=a\f{x(x-1)}2+ b\f{y(y-1)}2+c\f{z(z-1)}2+d\f{w(w-1)}2\Big\}\Big|.$$
The numbers $\f{x(x-1)}2\ (x\in\Bbb Z)$ are called triangular
numbers.
\par In 1828 Jacobi showed that
$$N(1,1,1,1;n)=8\sum_{d\mid n,4\nmid d}d.\tag 1.1$$
In 1847 Eisenstein (see [D]) gave formulas for $N(1,1,1,3;n)$ and
$N(1,1,1,5;n)$. From 1859 to 1866 Liouville made about 90
conjectures on $N(a,b,c,d;n)$ in a series of papers. Most
conjectures of Liouville have been proved. See [A1, A2,
AALW1-AALW4], Cooper's survey paper [C], Dickson's historical
comments [D] and Williams' book [W2].
\par
 Let
$$t'(a,b,c,d;n)=\Big|\Big\{(x,y,z,w)\in \Bbb N^4\ |\ n=a\f{x(x-1)}2+
b\f{y(y-1)}2+c\f{z(z-1)}2+d\f{w(w-1)}2\Big\}\Big|.$$ As $\f12
a(a+1)=\f12(-a-1)(-a-1+1)$ we have
$$t(a,b,c,d;n)=16t'(a,b,c,d;n).$$
In [L] Legendre stated that
$$t'(1,1,1,1;n)=\sigma(2n+1).\tag 1.2$$
In 2003, Williams [W1] showed that
$$t'(1,1,2,2;n)=\f 14\sum_{d\mid 4n+3}\big(d-(-1)^{\f{d-1}2}\big).$$
 For $a,b,c,d\in\Bbb N$ with $1<a+b+c+d\le 8$ let
$$C(a,b,c,d)=16+4i_1(i_1-1)i_2+8i_1i_3,$$
where $i_j$ is the number of elements in $\{a,b,c,d\}$ which are
equal to $j$. When $1<a+b+c+d\le 7$, in 2005 Adiga, Cooper and Han
[ACH] showed that
$$C(a,b,c,d)t'(a,b,c,d;n)=N(a,b,c,d;8n+a+b+c+d).\tag 1.3$$ When
$a+b+c+d=8$, in 2008 Baruah, Cooper and Han [BCH] proved that
$$C(a,b,c,d)t'(a,b,c,d;n)=N(a,b,c,d;8n+8)-N(a,b,c,d;2n+2).\tag 1.4$$
 In 2009,
Cooper [C] determined $t'(a,b,c,d;n)$ for $(a,b,c,d)=(1,1,1,3),\
(1,3,3,3),$ $(1,2,2,3),\ (1,3,6,6),\ (1,3,4,4),\ (1,1,2,6)$ and
$(1,3,12,12)$.
\par In this paper, by using some formulas for $N(a,b,c,d;n)$
in [A1, A2, AALW1-AALW4] and Ramanujan's theta functions we obtain
explicit
 formulas for $t(a,b,c,d;n)$ in the cases
 $(a,b,c,d)=(1,2,2,4),\ (1,2,4,4),\ (1,1,4,4),\ (1,4,4,4)$,
 $(1,3,3,9)$,
 $(1,1,9,9),\ (1,9,9,9)$, $(1,1,1,9)$, $(1,3,9,9)$ and $(1,1,3,9).$
\par For $m,n\in \Bbb N$ with $2\mid m$ and $2\nmid n$ we define
$$S_m(n) =\sum\Sb(r,s)\in \Bbb Z \times
\Bbb Z\\n=r^2+m{s^2}\\r\e 1\mod 4
\endSb r.$$
 As ${r^2}+2{s^2}\e 0,1,2,3,4,6\mod 8$ for $r,s\in\Bbb Z$, we see
that $S_2(n)=0\qtq{for} n\e 5,7\mod 8.$ Also, ${r^2}+4{s^2}\e
0,1\mod 4$ for $r,s\in\Bbb Z$ implies that $S_4(n)=0\qtq{for} n\e
2,3\mod 4.$ In this paper, following [AALW4] we also define
$$S(n)=\sum_{d\mid n} \f nd\Ls 2d,$$ where $\sls am$ is the
Legendre-Jacobi-Kronecker symbol.

\section*{2. Formulas for $t(1,3,9,9;n)$ and $t(1,1,3,9;n)$}
\par Ramanujan's theta functions $\varphi(q)$ and $\psi(q)$ are defined
by
$$\varphi(q)=\sum_{n=-\infty}^{\infty}q^{n^2}=1+2\sum_{n=1}^{\infty}
q^{n^2}\qtq{and} \psi(q)=\sum_{n=0}^{\infty}q^{n(n+1)/2}\ (|q|<1).$$
It is evident that for $|q|<1$,
$$\sum_{n=0}^{\infty}N(a,b,c,d;n)q^{n}=\varphi(q^a)
\varphi(q^b)\varphi(q^c)\varphi(q^d),$$
$$\sum_{n=0}^{\infty}t'(a,b,c,d;n)q^{n}=\psi(q^a)\psi(q^b)
\psi(q^c)\psi(q^d).$$ From [BCH, Lemma 4.1] we know that for
$|q|<1$,
 $$\varphi(q)=\varphi(q^4)+2q\psi(q^8)\tag 2.1$$
 and
 $$\psi(q)\psi(q^3)=\varphi(q^6)\psi(q^4)+q\psi(q^{12})
 \varphi(q^2).\tag
 2.2$$

\pro{Theorem 2.1} Let $n\in\Bbb N$. Then
$$N(1,3,9,9;8n+22)=40t'(1,3,9,9;n).$$
\endpro

Proof.  By (2.1), for $|q|<1$ we have
$$\varphi(q^k)=\varphi(q^{4k})+2q^k\psi(q^{8k})
=\varphi(q^{16k})+2q^{4k}\psi(q^{32k})+2q^k\psi(q^{8k}).$$ Thus, for
$|q|<1$ we have
$$\align&\sum_{n=0}^{\infty}N(1,3,9,9;n)q^{n}\\&=
\varphi(q)\varphi(q^3)\varphi(q^9)^2
\\&=(\varphi(q^{16})+2q^4\psi(q^{32})+2q\psi(q^8))(\varphi(q^{48})+2q^{12}\psi(q^{96})+
2q^3\psi(q^{24}))\\&\q\times(\varphi(q^{144})+2q^{36}\psi(q^{288})+
2q^9\psi(q^{72}))^2
\\&=\big(\varphi(q^{16})\varphi(q^{48})+4q^{16}\psi(q^{32})\psi(q^{96})
+2q\psi(q^{8})\varphi(q^{48})+2q^3 \psi(q^{24})\varphi(q^{16})
\\&\q+2q^4 \psi(q^{32})\varphi(q^{48})
+4q^4\psi(q^8)\psi(q^{24})
 +2q^{12}\psi(q^{96})\varphi(q^{16})+4q^7 \psi(q^{24})\psi(q^{32})
 \\&\q+4q^{13} \psi(q^8)\psi(q^{96})\big)
 \big(\varphi(q^{144})^2+4q^{72}\psi (q^{288})^2
+4q^{36}\varphi(q^{144})\psi(q^{288})+4q^{18}\psi(q^{72})^2
\\&\q+4q^9\varphi(q^{144})\psi(q^{72})
+8q^{45}\psi(q^{288})\psi(q^{72})\big).
\endalign$$
Since
$$\varphi(q^{8k})=1+2\sum_{n=1}^{\infty}
q^{8kn^2}\qtq{and} \psi(q^{8k})=\sum_{n=0}^{\infty}q^{8kn(n+1)/2}\
(|q|<1),$$ we see that for any nonnegative integers $k_1,k_2,m_1$
and $m_2$,
$$\varphi(q^{8k_1})^{m_1}\psi(q^{8k_2})^{m_2}
=\sum_{n=0}^{\infty}b_nq^{8n}\q (|q|<1).$$
 Now from the above we deduce that for $|q|<1$,

$$\align &\sum_{n=0}^{\infty}N(1,3,9,9;8n+6)q^{8n+6}
 \\&=2q\psi(q^{8})\varphi(q^{48})\cdot
 8q^{45}\psi(q^{288})\psi(q^{72})+2q^4 \psi(q^{32})\varphi(q^{48})
 \cdot 4q^{18}\psi(q^{72})^2
 \\&\q+4q^4\psi(q^8)\psi(q^{24})\cdot
 4q^{18}\psi(q^{72})^2+2q^{12}\psi(q^{96})\varphi(q^{16})
 \cdot 4q^{18}\psi(q^{72})^2
 \\&\q+4q^{13} \psi(q^8)\psi(q^{96})
 \cdot 4q^9\varphi(q^{144})\psi(q^{72})
 \\&=16q^{46}\varphi(q^{48})\psi(q^8)\psi(q^{72})\psi(q^{288})
+8q^{22}\varphi(q^{48})\psi(q^{32})\psi(q^{72})^2
+16q^{22}\psi(q^8)\psi(q^{24})\psi(q^{72})^2
\\&\q+8q^{30}\varphi(q^{16})\psi(q^{96})\psi(q^{72})^2
+16q^{22}\varphi(q^{144})\psi(q^8)\psi(q^{72})\psi(q^{96})
\endalign$$
and so
$$\align&\f 18\sum_{n=0}^{\infty}N(1,3,9,9;8n+6)q^{8n-16}
\\&=2q^{24}\varphi(q^{48})\psi(q^8)\psi(q^{72})\psi(q^{288})
+\varphi(q^{48})\psi(q^{32})\psi(q^{72})^2
+2\psi(q^8)\psi(q^{24})\psi(q^{72})^2
\\&\q+q^8\varphi(q^{16})\psi(q^{96})\psi(q^{72})^2
+2\varphi(q^{144})\psi(q^8)\psi(q^{72})\psi(q^{96}).
\endalign$$
Replacing $q$ with $q^{1/8}$ in the above we obtain
$$\align &\f 18\sum_{n=0}^{\infty}N(1,3,9,9;8n+22)q^n
\\&=\f 18\sum_{n=0}^{\infty}N(1,3,9,9;8n+6)q^{n-2}
\\&=2q^3\varphi(q^{6})\psi(q)\psi(q^{9})\psi(q^{36})
+\varphi(q^{6})\psi(q^{4})\psi(q^{9})^2
+2\psi(q)\psi(q^{3})\psi(q^{9})^2
\\&\q+q\varphi(q^{2})\psi(q^{12})\psi(q^{9})^2
+2\varphi(q^{18})\psi(q)\psi(q^{9})\psi(q^{12}).
\endalign$$
 Now applying (2.2) we get
$$\align&\f 18\sum_{n=0}^{\infty}N(1,3,9,9;8n+22)q^n
\\&=2\psi(q)\psi(q^3)\psi(q^9)^2
+\psi(q^9)^2\psi(q)\psi(q^3)+2\psi(q)\psi(q^9)\psi(q^3)\psi(q^9)
\\&=5\psi(q)\psi(q^3)\psi(q^9)^2
=5\sum_{n=0}^{\infty}t'(1,3,9,9;n)q^n.\endalign$$ Comparing the
coefficients of $q^n$ in the above expansion we obtain the result.
$\square$
\par\q\par For $n\in \Bbb N$ following [AALW3] we define
$$\aligned &A(n)=\sum_{d\mid n}d\Ls{12}{n/d},
\q B(n)=\sum_{d\mid n}d\Ls{-3}{d}\Ls{-4}{n/d},
\\&C(n)=\sum_{d\mid n}d\Ls{-3}{n/d}\Ls{-4}{d}
\qtq{and}D(n)=\sum_{d\mid n}d\Ls{12}{d}.\endaligned$$  Let $(a,b)$
be the greatest common divisor of integers $a$ and $b$. Suppose that
$n\in\Bbb N$ and $n=2^{\a}3^{\beta}n_1$, where $\a$ and $\beta$ are
non-negative integers, $n_1\in\Bbb N$ and
 $(n_1,6)=1$. From [AALW3, Theorem 3.1] we know that
 $$\aligned
&A(n)=2^{\a}3^{\beta}A(n_1), \q
B(n)=(-1)^{\a+\beta}2^{\a}\Ls{-3}{n_1}A(n_1),
\\&C(n)=(-1)^{\a+\beta+\f{n_1-1}2}3^{\beta}A(n_1)
\qtq{and} D(n)=\Ls 3{n_1}A(n_1).\endaligned\tag 2.3$$

 \pro{Lemma 2.1 ([A1, Theorem 1.2])} Let $n\in\Bbb N$. Then
$$\align N(1,3,9,9;n)=\cases 2A(n/3)+2B(n/3)-C(n/3)-D(n/3)
&\t{if $n\e0\mod3$,}
\\2A(n)-\f23B(n)+C(n)-\f13D(n)&\t{if $n\e1\mod3$,}
\\0&\t{if $n\e2\mod3$.}\endcases\endalign$$
\endpro

\pro{Theorem 2.2} Let $n\in\Bbb N$. Then
$$\align &t(1,3,9,9;n)\\&=\cases 0&\t{if $3\mid n-2$},
\\\f 43\sum_{d\mid 4n+11}d\sls 3d&\t{if $3\mid n$},
\\2\big(3^{\beta-1}\sls 3{n_1}-1\big)\sum_{d\mid n_1}d\sls 3d
&\t{if $3\mid n-1$
 and $4n+11=3^{\beta}n_1$ $(3\nmid n_1)$}.
\endcases\endalign$$
\endpro
Proof. By Theorem 2.1,
$$t(1,3,9,9;n)=16t'(1,3,9,9;n)=\f 25N(1,3,9,9;8n+22).$$
 Now applying Lemma 2.1 and (2.3) we deduce that
$$\align t(1,3,9,9;n)=\cases 0&\t{if $3\mid n-2$},
\\\f 43A(4n+11)&\t{if $3\mid n$},
\\2(3^{\beta-1}-\sls 3{n_1})A(n_1)&\t{if $3\mid n-1$
 and $4n+11=3^{\beta}n_1$ $(3\nmid n_1)$}.
\endcases\endalign$$
To see the result, we note that
$$A(m)=\sum_{d\mid m}d\ls{12}m\ls{12}d=\ls 3m\sum_{d\mid m}
d\ls 3d\qtq{for $m\in\Bbb N$ with $(6,m)=1$.}\square\tag 2.4$$
 \pro{Lemma 2.2 ([A1, Theorem 1.3])} Let $n\in\Bbb N$. Then
$$\align N(1,1,3,9;n)=\cases 2A(n/3)+2B(n/3)-C(n/3)-D(n/3)
&\t{if $n\e0\mod3$,}
\\4A(n)-\f43B(n)+2C(n)-\f23D(n)&\t{if $n\e1\mod3$,}
\\2A(n)-\f23B(n)+C(n)-\f13D(n)&\t{if $n\e2\mod3$.}
\endcases\endalign$$
\endpro
 \pro{Theorem 2.3} Let $n\in\Bbb
N$. Then
$$ t(1,1,3,9;n)=\cases -\f 83\sum_{d\mid 4n+7}d\sls 3d
&\t{if $3\mid n$},
\\\f 83\sum_{d\mid 4n+7}d\sls 3d  &\t{if $3\mid n-1$},
\\2\big(3^{\beta-1}\sls 3{n_1}-1\big)\sum_{d\mid n_1}d\sls 3d
&\t{if $3\mid n-2$ and $4n+7=3^{\beta}n_1$ $(3\nmid
n_1)$}.\endcases$$
\endpro
 Proof. Suppose $|q|<1$. Then clearly
 $$\sum_{n=0}^{\infty}N(1,1,3,9;n)q^{n}
 =\varphi(q)^2\varphi(q^3)\varphi(q^9).$$
Since $\varphi(q^k)=\varphi(q^{4k})+2q^k\psi(q^{8k})
=\varphi(q^{16k})+2q^{4k}\psi(q^{32k})+2q^k\psi(q^{8k})$ by (2.1),
we see that
$$\align&\varphi(q)^2\varphi(q^3)\varphi(q^9)
\\&=
\big(\varphi(q^{16})+2q^4\psi(q^{32})+2q\psi(q^8)\big)^2
\big(\varphi(q^{48})+2q^{12}\psi(q^{96})+2q^3\psi(q^{24})\big)
\\&\q\times
\big(\varphi(q^{144})+2q^{36}\psi(q^{288})+2q^9\psi(q^{72})\big)
\\&=\big(\varphi(q^{16})^2+4q^4\psi(q^{32})\varphi(q^{16})
+4q^8\psi(q^{32})^2
\\&\q+4q^2\psi(q^{8})^2+4q\varphi(q^{16})\psi(q^{8})+8q^{5}
\psi(q^{8})\psi(q^{32})\big)
\\&\q\times
\big(\varphi(q^{48})\varphi(q^{144})
+2q^{36}\varphi(q^{48})\psi(q^{288})
+2q^9\varphi(q^{48})\psi(q^{72})+2q^{12}\psi(q^{96})\varphi(q^{144})
\\&\q\ +4q^{48}\psi(q^{96})\psi(q^{288})+ 4q^{21}\psi(q^{96})\psi(q^{72})+2q^3\psi(q^{24})
\varphi(q^{144}\big)
\\&\q\ +4q^{39}\psi(q^{24})\psi(q^{288})
+4q^{12}\psi(q^{24})\psi(q^{72})\big).
\endalign$$
Note that $\varphi(q^{8k_1})^{m_1}\psi(q^{8k_2})^{m_2}
=\sum_{n=0}^{\infty}b_nq^{8n}$ for $|q|<1$ and any nonnegative
integers $k_1,k_2,m_1$ and $m_2$. From the above we deduce that
$$\align&\sum_{n=0}^{\infty}N(1,1,3,9;8n+6)q^{8n+6}
\\&=4q\varphi(q^{16})\psi(q^{8})\cdot 4q^{21}\psi(q^{96})\psi(q^{72})
+8q^{5} \psi(q^{8})\psi(q^{32})\cdot 2q^9\varphi(q^{48})\psi(q^{72})
\\&\q+4q^2\psi(q^{8})^2\cdot 2q^{36}\varphi(q^{48})\psi(q^{288})
+4q^2\psi(q^{8})^2\cdot 2q^{12}\psi(q^{96})\varphi(q^{144})
\\&\q +4q^2\psi(q^{8})^2\cdot 4q^{12}\psi(q^{24})\psi(q^{72})
\endalign$$ and so
$$\align&\f 18\sum_{n=0}^{\infty}N(1,1,3,9;8n+6)q^{8n-8}
\\&=2q^{8}\varphi(q^{16})\psi(q^8)\psi(q^{72})\psi(q^{96})
+2\psi(q^8)\psi(q^{32})\varphi(q^{48})\psi(q^{72})
\\&\q+q^{24}\psi(q^8)^2 \varphi(q^{48})\psi(q^{288})
+\psi(q^8)^2\psi(q^{96})\varphi(q^{144})
+2\psi(q^8)^2\psi(q^{24})\psi(q^{72}).\endalign$$
 Replacing $q$ with
$q^{1/8}$ in the above we obtain
$$\align
&\f 18\sum_{n=0}^{\infty}N(1,1,3,9;8n+14)q^n
\\&=\f 18\sum_{n=0}^{\infty}N(1,1,3,9;8n+6)q^{n-1}
\\&=2q\varphi(q^{2})\psi(q)\psi(q^{9})\psi(q^{12})
+2\psi(q)\psi(q^{4})\varphi(q^{6})\psi(q^{9})+q^3\psi(q)^2
\varphi(q^{6})\psi(q^{36})
\\&\q+\psi(q)^2\psi(q^{12})\varphi(q^{18})
+2\psi(q)^2\psi(q^{3})\psi(q^{9}).
\endalign$$
Now applying (2.2) we get
$$\align&\f 18\sum_{n=0}^{\infty}N(1,1,3,9;8n+14)q^n
\\&=2\psi(q)^2\psi(q^3)\psi(q^9)
+\psi(q)^2\psi(q^3)\psi(q^9)+2\psi(q)^2\psi(q^3)\psi(q^9)
\\&=5\psi(q)^2\psi(q^3)\psi(q^9)
=5\sum_{n=0}^{\infty}t'(1,1,3,9;n)q^n =\f
5{16}\sum_{n=0}^{\infty}t(1,1,3,9;n)q^n.\endalign$$ Comparing the
coefficients of $q^n$ we obtain
$$t(1,1,3,9;n)=\f 25N(1,1,3,9;8n+14).$$
Now applying Lemma 2.2, (2.3) and (2.4) we deduce the result.
$\square$

\section*{3. Formulas for $t(1,1,4,4;n)$, $t(1,4,4,4;n)$, $t(1,2,2,4;n)$
and $t(1,2,4,4;n)$}
 \pro{Lemma 3.1}
Let $a,b,c,d,n\in\Bbb N$. Then
$$\align &t(a,b,c,d;n)\\&=N(a,b,c,d;8n+a+b+c+d)-N(a,b,c,4d;8n+a+b+c+d)
\\&\q-N(a,b,4c,d;8n+a+b+c+d)+N(a,b,4c,4d;8n+a+b+c+d)
\\&\q-N(a,4b,c,d;8n+a+b+c+d)+N(a,4b,c,4d;8n+a+b+c+d)
\\&\q+N(a,4b,4c,d;8n+a+b+c+d)-N(a,4b,4c,4d;8n+a+b+c+d)
\\&\q-N(4a,b,c,d;8n+a+b+c+d)+N(4a,b,c,4d;8n+a+b+c+d)
\\&\q+N(4a,b,4c,d;8n+a+b+c+d)-N(4a,b,4c,4d;8n+a+b+c+d)
\\&\q+N(4a,4b,c,d;8n+a+b+c+d)-N(4a,4b,c,4d;8n+a+b+c+d)
\\&\q-N(4a,4b,4c,d;8n+a+b+c+d)+N(4a,4b,4c,4d;8n+a+b+c+d).\endalign$$
\endpro
Proof. It is clear that
$$\align &t(a,b,c,d;n)\\&=\big|\{(x,y,z,w)\in
\Bbb Z^4\bigm|n=a\f{x(x-1)}2+
b\f{y(y-1)}2+c\f{z(z-1)}2+d\f{w(w-1)}2\}\big|
\\&=\big|\{(x,y,z,w)\in \Bbb Z^4
\bigm|8n+a+b+c+d\\&\qq\qq=a{(2x-1)^2}+b{(2y-1)^2}+c{(2z-1)^2}
+d{(2w-1)^2}\}\big|
\\&=\big|\{(x,y,z,w)\in \Bbb Z^4\bigm|n=a{x^2}+b{y^2}+c{z^2}+d{w^2},
2\mid{xyzw-1}\}\big|
\\&=\big|\{(x,y,z,w)\in \Bbb Z^4\bigm|n=a{x^2}+b{y^2}+c{z^2}+d{w^2},
2\mid{yzw-1}\}\big|
\\&\q-\big|\{(x,y,z,w)\in \Bbb Z^4\bigm|n=4a{x^2}+b{y^2}+c{z^2}+d{w^2},
2\mid{yzw-1}\}\big|
\\&=\big|\{(x,y,z,w)\in \Bbb Z^4\bigm|n=a{x^2}+b{y^2}+c{z^2}+d{w^2},
2\mid{zw-1}\}\big|
\\&\q-\big|\{(x,y,z,w)\in \Bbb Z^4\bigm|n=a{x^2}+4b{y^2}+c{z^2}+d{w^2},
2\mid{zw-1}\}\big|
\\&\q-\big|\{(x,y,z,w)\in \Bbb Z^4\bigm|n=4a{x^2}+b{y^2}+c{z^2}+d{w^2},
2\mid{zw-1}\}\big|
\\&\q+\big|\{(x,y,z,w)\in \Bbb Z^4\bigm|n=4a{x^2}+4b{y^2}+c{z^2}+d{w^2},
2\mid{zw-1}\}\big|
\\&=\big|\{(x,y,z,w)\in \Bbb Z^4\bigm|n=a{x^2}+b{y^2}+c{z^2}+d{w^2},
2\mid{w-1}\}\big|
\\&\q-\big|\{(x,y,z,w)\in \Bbb Z^4\bigm|n=a{x^2}+b{y^2}+4c{z^2}+d{w^2},
2\mid{w-1}\}\big|
\\&\q-\big|\{(x,y,z,w)\in \Bbb Z^4\bigm|n=a{x^2}+4b{y^2}+c{z^2}+d{w^2},
2\mid{w-1}\}\big|
\\&\q+\big|\{(x,y,z,w)\in \Bbb Z^4\bigm|n=a{x^2}+4b{y^2}+4c{z^2}+d{w^2},
2\mid{w-1}\}\big|
\\&\q-\big|\{(x,y,z,w)\in \Bbb Z^4\bigm|n=4a{x^2}+b{y^2}+c{z^2}+d{w^2},
2\mid{w-1}\}\big|
\\&\q+\big|\{(x,y,z,w)\in \Bbb Z^4\bigm|n=4a{x^2}+b{y^2}+4c{z^2}+d{w^2},
2\mid{w-1}\}\big|
\\&\q+\big|\{(x,y,z,w)\in \Bbb Z^4\bigm|n=4a{x^2}+4b{y^2}+c{z^2}+d{w^2},
2\mid{w-1}\}\big|
\\&\q-\big|\{(x,y,z,w)\in \Bbb Z^4\bigm|n=4a{x^2}+4b{y^2}+4c{z^2}+d{w^2},
2\mid{w-1}\}\big|
\\&=\big|\{(x,y,z,w)\in \Bbb Z^4\bigm|n=a{x^2}+b{y^2}+c{z^2}+d{w^2}\}\big|
\\&\q-\big|\{(x,y,z,w)\in \Bbb Z^4\bigm|n=a{x^2}+b{y^2}+c{z^2}+4d{w^2}\}\big|
\\&\q-\big|\{(x,y,z,w)\in \Bbb Z^4\bigm|n=a{x^2}+b{y^2}+4c{z^2}+d{w^2}\}\big|
\\&\q+\big|\{(x,y,z,w)\in \Bbb Z^4\bigm|n=a{x^2}+b{y^2}+4c{z^2}+4d{w^2}\}\big|
\\&\q-\big|\{(x,y,z,w)\in \Bbb Z^4\bigm|n=a{x^2}+4b{y^2}+c{z^2}+d{w^2}\}\big|
\\&\q+\big|\{(x,y,z,w)\in \Bbb Z^4\bigm|n=a{x^2}+4b{y^2}+c{z^2}+4d{w^2}\}\big|
\\&\q+\big|\{(x,y,z,w)\in \Bbb Z^4\bigm|n=a{x^2}+4b{y^2}+4c{z^2}+d{w^2}\}\big|
\\&\q-\big|\{(x,y,z,w)\in \Bbb Z^4\bigm|n=a{x^2}+4b{y^2}+4c{z^2}+4d{w^2}\}\big|
\\&\q-\big|\{(x,y,z,w)\in \Bbb Z^4\bigm|n=4a{x^2}+b{y^2}+c{z^2}+d{w^2}\}\big|
\\&\q+\big|\{(x,y,z,w)\in \Bbb Z^4\bigm|n=4a{x^2}+b{y^2}+c{z^2}+4d{w^2}\}\big|
\\&\q+\big|\{(x,y,z,w)\in \Bbb Z^4\bigm|n=4a{x^2}+b{y^2}+4c{z^2}+d{w^2}\}\big|
\\&\q-\big|\{(x,y,z,w)\in \Bbb Z^4\bigm|n=4a{x^2}+b{y^2}+4c{z^2}+4d{w^2}\}\big|
\\&\q+\big|\{(x,y,z,w)\in \Bbb Z^4\bigm|n=4a{x^2}+4b{y^2}+c{z^2}+d{w^2}\}\big|
\\&\q-\big|\{(x,y,z,w)\in \Bbb Z^4\bigm|n=4a{x^2}+4b{y^2}+c{z^2}+4d{w^2}\}\big|
\\&\q-\big|\{(x,y,z,w)\in \Bbb Z^4\bigm|n=4a{x^2}+4b{y^2}+4c{z^2}+d{w^2}\}\big|
\\&\q+\big|\{(x,y,z,w)\in \Bbb Z^4\bigm|
n=4a{x^2}+4b{y^2}+4c{z^2}+4d{w^2}\}\big|.
\endalign$$
Thus the result follows. $\square$
\par For general positive integer $n$, in a series of papers A. Alaca, S. Alaca, M.F. Lemire and
K.S. Williams (see [AALW1, AALW2, AALW4]) gave explicit formulas for
$N(a,b,c,d;n)$ in the cases $(a,b,c,d)=(1,1,4,4)$, $(1,1,16,16)$,
$(1,1,4,16)$, $(1,4,4,4), (1,4,16,16),$ $(1,4,4,$ $16),\ (1,2,2,4),\
(1,2,2,16)$, $(1,2,16,16),\ (1,2,4,16),\ (1,2,4,8), \ (1,2,4,4),\
(1,2,8,16),$ $(1,4,8,8)$ and $(1,8,8,16)$. Based on Lemma 3.1, we
need some special results in [AALW1, AALW2, AALW4] to prove our
formulas for $t(1,1,4,4;n)$, $t(1,4,4,4;n)$, $t(1,2,2,4;n)$ and
$t(1,2,4,4;n)$.
 \pro{Lemma 3.2 ([AALW1,
Theorem 1.11])} Let $n\in\Bbb N$ with $n\e 2\mod 4$. Then
$N(1,1,4,4;n)=4\sigma(n/2).$
\endpro

\pro{Lemma 3.3 ([AALW2, Theorems 4.6 and 4.8])} Let $n\in\Bbb N$ and
$n\e 2\mod 8$. Then
$$N(1,1,16,16;n)=N(1,1,4,16;n)
=2\sigma\big(\f n2\big)+2\ls 2{n/2}S_4(\f n2).$$
\endpro

\pro{Theorem 3.1} Let $n\in\Bbb N$. Then
$$t(1,1,4,4;n)=2\big(\sigma(4n+5)+(-1)^nS_4(4n+5)\big).$$
\endpro
Proof. Since $x^2\not\e 2\mod 4$ for $x\in\Bbb Z$, from Lemma 3.1 we
see that
$$\aligned &t(1,1,4,4;n)\\&=N(1,1,4,4;8n+10)-N(1,1,4,16;8n+10)
\\&\q-N(1,1,16,4;8n+10)+N(1,1,16,16;8n+10)-N(1,4,4,4;8n+10)
\\&\q+N(1,4,4,16;8n+10)+N(1,4,16,4;8n+10)-N(1,4,16,16;8n+10)
\\&\q-N(4,1,4,4;8n+10)+N(4,1,4,16;8n+10)+N(4,1,16,4;8n+10)
\\&\q-N(4,1,16,16;8n+10)+N(4,16,4,4;8n+10)-N(4,4,4,16;8n+10)
\\&\q-N(4,4,16,4;8n+10)+N(4,4,16,16;8n+10)
\\&=N(1,1,4,4;8n+10)-2N(1,1,4,16;8n+10)+N(1,1,16,16;8n+10).
\endaligned$$
Now applying Lemmas 3.2 and 3.3 we obtain
$$\align
t(1,1,4,4;n)&=4\sigma(4n+5)-2\Big(2\sigma(4n+5)+2\Ls2{4n+5}S_4(4n+5)\Big)
\\&\q+2\sigma(4n+5)+2\Ls2{4n+5}S_4(4n+5)
\\&=2\Big(\sigma(4n+5)-\Ls2{4n+5}S_4(4n+5)\Big).\endalign$$
This yields the result. $\square$

\pro{Lemma 3.4 ([AALW1, Theorem 1.18])} Let $n\in\Bbb N$ and $n\e
1\mod 4$. Then
$$N(1,4,4,4;n)= 2\sigma(n).$$
\endpro
\pro{Lemma 3.5 ([AALW2, Theorem 4.5])} Let $n\in\Bbb N$ and $n\e
1\mod 4$. Then
$$N(1,4,16,16;n)=\f12\sigma(n)+
\f12(2+(-1)^{\f{n-1}4})S_4(n).$$
\endpro

\pro{Lemma 3.6 ([AALW2, Theorem 4.7])} Let $n\in\Bbb N$ and $n\e
1\mod 4$. Then
$$N(1,4,4,16;n)=\sigma(n)+S_4(n).$$
\endpro
\pro{Theorem 3.2} Let $n\in\Bbb N$. Then
$$ t(1,4,4,4;n)=\f 12\Big(\sigma(8n+13)-3S_4(8n+13)\Big).$$
\endpro
Proof. Since $x^2\e 0,1,4\mod 8$ for $x\in\Bbb Z$, using Lemma 3.1
we see that
$$\aligned &t(1,4,4,4;n)\\&=N(1,4,4,4;8n+13)-N(1,4,4,16;8n+13)
\\&\q-N(1,4,16,4;8n+13)+N(1,4,16,16;8n+13)-N(1,16,4,4;8n+13)
\\&\q+N(1,16,4,16;8n+13)+N(1,16,16,4;8n+13)-N(1,16,16,16;8n+13)
\\&\q-N(4,4,4,4;8n+13)+N(4,4,4,16;8n+13)+N(4,4,16,4;8n+13)
\\&\q-N(4,4,16,16;8n+13)+N(4,16,4,4;8n+13)-N(4,16,4,16;8n+13)
\\&\q-N(4,16,16,4;8n+13)+N(4,16,16,16;8n+13)
\\&=N(1,4,4,4;8n+13)-3N(1,4,4,16;8n+13)+3N(1,4,16,16;8n+13).
\endaligned$$
Now applying Lemmas 3.4, 3.5 and 3.6 we obtain
$$\align &t(1,4,4,4;n)\\&=2\sigma(8n+13)-3(\sigma(8n+13)+S_4(8n+13))
+\f 32(\sigma(8n+13)+S_4(8n+13))
\\&=\f12\Big(\sigma(8n+13)-3S_4(8n+13)\Big).\endalign$$
This proves the theorem. $\square$

\pro{Lemma 3.7 ([AALW1, Theorem 1.14])} Let $n\in\Bbb N$ with
$2\nmid n$. Then $$N(1,2,2,4;n)=2\sigma(n).$$
\endpro

 \pro{Lemma 3.8 ([AALW2, Theorems 4.9, 4.11 and
4.13])} Let $n\in\Bbb N$ and $n\e 1\mod 8$. Then
$$N(1,2,2,16;n)=N(1,8,8,16;n)=N(1,2,8,16;n)=\sigma(n)+S_4(n).$$
\endpro
\pro{Lemma 3.9 ([AALW2, Theorems 4.1 and 4.4])} Let $n\in\Bbb N$ and
$n\e 1\mod 4$. Then
$$N(1,2,4,8;n)=N(1,4,8,8;n)=\sigma(n)+(-1)^{\f{n-1}4}S_4(n).$$
\endpro

\pro{Theorem 3.3} Let $n\in\Bbb N$. Then
$$t(1,2,2,4;n)=\sigma(8n+9)-S_4(8n+9).$$\endpro
Proof. From Lemma 3.1 we have
$$\aligned &t(1,2,2,4;n)\\&=N(1,2,2,4;8n+9)-N(1,2,2,16;8n+9)
\\&\q-N(1,2,4,8;8n+9)+N(1,2,8,16;8n+9)-N(1,2,4,8;8n+9)
\\&\q+N(1,2,8,16;8n+9)+N(1,4,8,8;8n+9)-N(1,8,8,16;8n+9)
\\&\q-N(4,2,2,4;8n+9)+N(4,2,2,16;8n+9)+N(4,2,8,4;8n+9)
\\&\q-N(4,2,8,8;8n+9)+N(4,8,4,4;8n+9)-N(4,8,2,16;8n+9)
\\&\q-N(4,8,8,4;8n+9)+N(4,8,8,16;8n+9)
\\&=N(1,2,2,4;8n+9)-N(1,2,2,16;8n+9)-2N(1,2,4,8;8n+9)
\\&\q+2N(1,2,8,16;8n+9)+N(1,4,8,8;8n+9)-N(1,8,8,16;8n+9).
\endaligned$$ Now applying
Lemmas 3.7, 3.8 and 3.9 we obtain
$$\align &t(1,2,2,4;n)
\\&=2\sigma(8n+9)-(\sigma(8n+9)+S_4(8n+9))-2(\sigma(8n+9)+S_4(8n+9))
\\&\q+2(\sigma(8n+9)+S_4(8n+9))+\sigma(8n+9)+S_4(8n+9)-(\sigma(8n+9)+S_4(8n+9))
\\&=\sigma(8n+9)-S_4(8n+9),
\endalign$$
which completes the proof. $\square$
\endpro
\pro{Lemma 3.10 ([AALW2, Theorems 4.17 and 4.18])} Let $n\in\Bbb N$
and $n\e 1,3\mod8$. Then
$$N(1,2,4,16;n)=N(1,2,16,16;n)=S(n)+S_2(n).$$
\endpro

\pro{Lemma 3.11 ([AALW4, Theorem 5.4])} Let $n\in\Bbb N$ with
$2\nmid n$. Then
$$N(1,2,4,4;n)=
 2S(n).$$
\endpro
\pro{Theorem 3.4} Let $n\in\Bbb N$. Then
$$t(1,2,4,4;n)=-\sum_{d\mid 8n+11}d\Ls 2d-S_2(8n+11).$$
\endpro
Proof. Since $x^2\e 0,1\mod 4$ for $x\in\Bbb Z$, from Lemma 3.1 we
see that
$$\align &t(1,2,4,4;n)\\&=N(1,2,4,4;8n+11)-N(1,2,4,16;8n+11)
\\&\q-N(1,2,16,4;8n+11)+N(1,2,16,16;8n+11)-N(1,8,4,4;8n+11)
\\&\q+N(1,8,4,16;8n+11)+N(1,8,16,4;8n+11)-N(1,8,16,16;8n+11)
\\&\q-N(4,2,4,4;8n+11)+N(4,2,4,16;8n+11)+N(4,2,16,4;8n+11)
\\&\q-N(4,2,16,16;8n+11)+N(4,8,4,4;8n+11)-N(4,8,4,16;8n+11)
\\&\q-N(4,8,16,4;8n+11)+N(4,8,16,16;8n+11)
\\&=N(1,2,4,4;8n+11)-2N(1,2,4,16;8n+11)+N(1,2,16,16;8n+11).
\endalign$$
Now applying Lemmas 3.10 and 3.11 we obtain $$\align
t(1,2,4,4;n)&=2S(8n+11)-2(S(8n+11)+S_2(8n+11))+(S(8n+11)+S_2(8n+11))
\\&=S(8n+11)-S_2(8n+11).
\endalign$$
Since
$$S(8n+11)=\sum_{d\mid 8n+11}\f nd\Ls 2d=\sum_{d\mid 8n+11}d\Ls
2{(8n+11)/d}=-\sum_{d\mid 8n+11}d\Ls 2d,$$ from the above we
 deduce the result. $\square$

 \section*{4. Formulas for $t(1,3,3,9;n)$, $t(1,1,9,9;n)$, $t(1,9,9,9;n)$ and $t(1,1,1,9;n)$}
\par For $a,b,c,d,n\in\Bbb N$ let
$$N_0(a,b,c,d;n)=\big|\big\{(x,y,z,w)\in\Bbb Z^4\bigm|
n=ax^2+by^2+cz^2+dw^2,\ 2\nmid xyzw\big\}\big|.$$ From the proof of
Lemma 3.1 we know that
$$t(a,b,c,d;n)=N_0(a,b,c,d;8n+a+b+c+d).\tag 4.1$$
\pro{Lemma 4.1} Let $n\in\Bbb N$ and $n+1=2^{\alpha}3^{\beta}n_1$
with $(6,n_1)=1$. Then
$$t(1,1,3,3;n)=\cases 4\sigma(n_1)&\t{if $2\mid n$,}
\\2^{\a+4}\sigma(n_1)&\t{if $2\nmid n$.}
\endcases$$\endpro
Proof. By [BCH, Theorem 1.5],
$$t(1,1,3,3;n)=16t'(1,1,3,3;n)=\cases 4N(1,1,3,3;n+1)
&\t{if $2\mid n$,}
\\2(N(1,1,3,3;2n+2)-N(1,1,3,3;n+1))&\t{if $2\nmid n$.}
\endcases$$
 Ramanujan (see [Be, pp. 114,223]) gave theta function identities that yields
  the following Liouville's conjecture (see [D]):
$$N(1,1,3,3;n+1)=\cases 16\sigma(n_1)&\t{if $2\mid n$,}
\\4(2^{\a+1}-3)\sigma(n_1)&\t{if $2\nmid n$.}
\endcases$$
Since $2n+2=2^{\a+1}3^{\beta}n_1$, combining all the above yields
the result. $\square$

 \pro{Theorem 4.1} Let $n\in\Bbb N$ and $n+2=2^{\a}3^{\beta}n_1$ with
  $(6,n_1)=1$. Then
$$t(1,3,3,9;n)=\cases 0&\t{if $n\e 2,5\mod 6$,}
\\ 16\sigma(n_1)&\t{if $n\e 1\mod 6$,}
\\2^{\a+4}\sigma(n_1)&\t{if $n\e 4\mod 6$,}
\\8\sigma(n_1)&\t{if $n\e 3\mod 6$,}
\\2^{\a+3}\sigma(n_1)&\t{if $n\e 0\mod 6$.}
\endcases$$
\endpro
Proof. From (4.1) we know that $t(1,3,3,9;n)=N_0(1,3,3,9;8n+16)$. If
$3\mid n-2$, then $8n+16\e 2\mod 3$. Since $x^2\not\e 2\mod 3$ for
any $x\in\Bbb Z$, we get $t(1,3,3,9;n)=N_0(1,3,3,9;8n+16)=0$. If
$3\mid n-1$, then $3\mid 8n+16$ and so
$$\align t(1,3,3,9;n)
&=N_0(1,3,3,9;8n+16) \\&=\big|\big\{(x,y,z,w)\in\Bbb Z^4\bigm|
8n+16=(3x)^2+3y^2+3z^2+9w^2,\ 2\nmid xyzw\big\}\big|
\\&=\big|\big\{(x,y,z,w)\in\Bbb Z^4\bigm|
\f{8n+16}3=3x^2+y^2+z^2+3w^2,\ 2\nmid xyzw\big\}\big|
\\&=N_0(1,1,3,3;8(n-1)/3+8)=t(1,1,3,3;(n-1)/3).
\endalign$$
 If $3\mid n$, since $x^2+y^2\e 8n+16\e 1\mod 3$ implies $3\mid
 x$ or $3\mid y$ we see that
 $$\align t(1,1,3,3;n+1)
 &=N_0(1,1,3,3;8n+16)
 \\&=\big|\big\{(x,y,z,w)\in\Bbb Z^4\bigm|
8n+16=(3x)^2+y^2+3z^2+3w^2,\ 2\nmid xyzw\big\}\big|
\\&\q+\big|\big\{(x,y,z,w)\in\Bbb Z^4\bigm|
8n+16=x^2+(3y)^2+3z^2+3w^2,\ 2\nmid xyzw\big\}\big|
\\&=2N_0(1,3,3,9;8n+16)=2t(1,3,3,9;n).\endalign$$
Now combining the above with Lemma 4.1 yields the result. $\square$
\par\q\par
 For $k\in\Bbb N$
and variable $q$ with $|q|<1$ define
$$\align &E_k=E_k(q):=\prod_{n=1}^{\infty}(1-q^{kn}),
\\& qE_6^4=\sum_{n=1}^{\infty}c(n)q^n\qtq{and}
\f{E_2^{17}E_3}{E_1^7E_4^6E_6}=\sum_{n=0}^{\infty}a(n)q^n.
\endalign$$ From [KF, p.374] or [M, p.121] we know that
$$c(n)=\f 13\sum\Sb x,y\in\Bbb Z
\\ n=x^2+3xy+3y^2\\ 3\mid x-2,\ 2\mid y-1\endSb
(-1)^xx.$$ Thus,
$$c(n)=\f 13\sum\Sb x,y\in\Bbb Z\\n=x^2+3x(1+2y)+3(1+2y)^2
\\x\e 2\mod 3\endSb(-1)^xx
=\f 13\sum\Sb x,y\in\Bbb Z\\4n=x^2+3(x+2+4y)^2
\\x\e 2\mod 3\endSb(-1)^xx$$ and so
$$c(n)=\f 13\sum\Sb 4n=a^2+3b^2 \ (a,b\in\Bbb Z)
\\a\e 2\mod 3,b\e a+2\mod 4\endSb(-1)^aa.\tag 4.2$$
 \pro{Lemma
4.2 ([A1, Theorems 1.5 and 1.6])} For $n\in\Bbb N$ we have
$$ N(1,1,9,9;n)=\cases 4\sigma(n)-8\sigma(
n/2)&\t{if $n\e2,4\mod 6$,}\\\f 43\sigma(n)&\t{if $n\e5\mod
6$,}\\8\sigma(n/9)-32\sigma( n/36)&\t{if $n\e0\mod 6$}\endcases$$
and
$$ N(1,9,9,9;n)=\cases 8\sigma(n/9)&\t{if $n\e3\mod 6$,}\\2\sigma(n)-4\sigma(
n/2)&\t{if $n\e4\mod 6$,}\\8\sigma(n/9)-32\sigma(n/36)&\t{if
$n\e0\mod 6$.}\endcases$$
\endpro
\pro{Lemma 4.3 ([A2, Theorems 2.5 and 2.10])} For $n\in\Bbb N$ with
$4\mid n$ we have
$$
N(1,1,36,36;n)=\cases \f43\sigma(n/4)-\f{16}3\sigma(n/16)+\f
83c(n/4)&\t{if $n\e4\mod {12}$,}
\\\f43\sigma(n/4)-\f{16}3\sigma(n/16)&\t{if $n\e8\mod {12}$,}
\\8\sigma(n/36)-32\sigma(n/144)&\t{if $n\e0\mod {12}$}
\endcases$$
and
$$ N(1,4,36,36;n)=\cases \f43\sigma(n/4)-\f{16}3\sigma(n/16)+\f 83c(n/4)&\t{if
$n\e4\mod {12}$,}
\\\f43\sigma(n/4)-\f{16}3\sigma(n/16)&\t{if $n\e8\mod{12}$,}
\\8\sigma(n/36)-32\sigma(n/144)&\t{if $n\e0\mod {12}$.}
\endcases$$\endpro
\pro{Lemma 4.4 ([A2, Theorem 2.4])} For $n\in\Bbb N$  with $4\mid n$
we have
$$ N(1,1,9,36;n)=\cases \f43\sigma(n/4)-\f{16}3\sigma(n/16)+\f 83c(n/4)&\t{if $n\e4\mod
{12}$,}
\\\f43\sigma(n/4)-\f{16}3\sigma(n/16)&\t{if $n\e8\mod {12}$,}
\\8\sigma(n/36)-32\sigma(n/144)&\t{if $n\e0\mod {12}$.}
\endcases$$
\endpro
\pro{Lemma 4.5 ([A2, Theorem 2.8])} For $n\in\Bbb N$  with $4\mid n$
we have
$$ N(1,4,9,9;n)=\cases \f43\sigma(n/4)-\f{16}3\sigma(n/16)+\f 83c(n/4)&\t{if $n\e4\mod
{12}$,}
\\\f43\sigma(n/4)-\f{16}3\sigma(n/16)&\t{if $n\e8\mod {12}$,}
\\8\sigma(n/36)-32\sigma(n/144)&\t{if $n\e0\mod {12}$.}
\endcases$$\endpro
\pro{Lemma 4.6 ([A2, Theorem 2.9])} For $n\in\Bbb N$  with $4\mid n$
we have
$$ N(1,4,9,36;n)=\cases \f43\sigma(n/4)-\f{16}3\sigma(n/16)+\f 83c(n/4)
&\t{if $n\e4\mod{12}$,}
\\\f43\sigma(n/4)-\f{16}3\sigma(n/16)&\t{if $n\e8\mod {12}$,}
\\8\sigma(n/36)-32\sigma(n/144)&\t{if $n\e0\mod {12}$.}
\endcases$$\endpro
\pro{Lemma 4.7 ([A2, Theorem 2.15])} For $n\in\Bbb N$  with $4\mid
n$ we have
$$ N(4,4,9,9;n)=\cases
\f43\sigma(n/4)-\f{16}3\sigma(n/16)+\f83c(n/4)&\t{if $n\e4\mod
{12}$,}
\\\f43\sigma(n/4)-\f{16}3\sigma(n/16)&\t{if $n\e8\mod{12}$,}
\\8\sigma(n/36)-32\sigma(n/144)&\t{if $n\e0\mod {12}$.}\endcases$$

\pro{Lemma 4.8 ([A2, Theorem 2.16])} For $n\in\Bbb N$  with $4\mid
n$ we have
 $$N(4,4,9,36;n)=\cases
\f43\sigma(n/4)-\f{16}3\sigma(n/16)+\f83c(n/4)&\t{if $n\e4\mod
{12}$,}
\\\f43\sigma(n/4)-\f{16}3\sigma(n/16)&\t{if $n\e8\mod {12}$,}
\\8\sigma(n/36)-32\sigma(n/144)&\t{if $n\e0\mod {12}$.}
\endcases$$
\endpro
\pro{Theorem 4.2} Let $n\in\Bbb N.$ Then
$$\aligned t(1,1,9,9;n)=\cases \f83\sigma(2n+5)&\t{if $n\e0\mod 3$,}
\\{16}\sigma(\f{2n+5}9)&\t{if $n\e2\mod 9$,}
\\0&\t{if $n\e5,8\mod 9$,}
\\\f 83(\sigma(2n+5)-c(2n+5))&\t{if $n\e1\mod 3$.}
\endcases\endaligned$$
\endpro
Proof. For $n\e 2\mod 3$ we see that $3\mid 8n+20$ and so
$$\align &t(1,1,9,9,;n)\\&=N_0(1,1,9,9;8n+20)
\\&=\big|\big\{(x,y,z,w)\in\Bbb Z^4\bigm|
8n+20=(3x)^2+(3y)^2+9z^2+9w^2,\ 2\nmid xyzw\big\}\big|
\\&=\cases 0&\t{if $9\nmid n-2$,}
\\N_0(1,1,1,1;\f{8n+20}9)=t(1,1,1,1;\f{n-2}9)=16\sigma(\f{2n+5}9)
&\t{if $9\mid n-2$.}
\endcases\endalign$$
\par Now assume $n\e 0,1\mod 3$.
By Lemma 3.1,
$$\align &t(1,1,9,9;n)\\&=N(1,1,9,9;8n+20)-2N(1,1,9,36;8n+20)+N(1,1,36,36;8n+20)
\\&\q-2N(1,4,9,9;8n+20)+4N(1,4,9,36;8n+20)-2N(1,4,36,36;8n+20)
\\&\q+N(4,4,9,9;8n+20)-2N(4,4,9,36;8n+20)+N(1,1,9,9;2n+5).\endalign$$
 For $n\e0\mod3$ applying Lemmas 4.2-4.8 we see that
$$\align
&t(1,1,9,9;n)\\&=4\sigma(8n+20)-8\sigma(\f{8n+20}2)-2(\f43\sigma(\f{8n+20}4)
-\f{16}3\sigma(\f{8n+20}{16}))
\\&\q+\f43\sigma(\f{8n+20}4)-\f{16}3\sigma(\f{8n+20}{16})-2(\f43\sigma(\f{8n+20}4)
-\f{16}3\sigma(\f{8n+20}{16}))
\\&\q+4(\f43\sigma(\f{8n+20}4)-\f{16}3\sigma(\f{8n+20}{16}))-2(\f43\sigma(\f{8n+20}4)
-\f{16}3\sigma(\f{8n+20}{16}))
\\&\q+\f43\sigma(\f{8n+20}4)-\f{16}3\sigma(\f{8n+20}{16})-2(\f43\sigma(\f{8n+20}4)
-\f{16}3\sigma(\f{8n+20}{16}))+\f43\sigma(2n+5)
\\&=28\sigma(2n+5)-24\sigma(2n+5)-\f83\sigma(2n+5)+\f{32}3\sigma(\f{4n+5}4)+\f43\sigma(2n+5)
\\&=\f83\sigma(2n+5).\endalign$$
For $n\e 1\mod 3$, applying Lemmas 4.2-4.8 we find that
$$\align &t(1,1,9,9;n)\\&=4\sigma(8n+20)-8\sigma(\f{8n+20}2)-2(\f43\sigma(\f{8n+20}{4})
-\f{16}3\sigma(\f{8n+20}{16})+\f83c(\f{8n+20}4))
\\&\q+\f43\sigma(\f{8n+20}{4})
-\f{16}3\sigma(\f{8n+20}{16})+\f83c(\f{8n+20}4)-2(\f43\sigma(\f{8n+20}{4})
-\f{16}3\sigma(\f{8n+20}{16})
\\&\q+\f83c(\f{8n+20}4))
+4(\f43\sigma(\f{8n+20}{4})
-\f{16}3\sigma(\f{8n+20}{16})+\f83c(\f{8n+20}4))-2(\f43\sigma(\f{8n+20}{4})
\\&\q-\f{16}3\sigma(\f{8n+20}{16})+\f83c(\f{8n+20}4))+\f43\sigma(\f{8n+20}{4})
-\f{16}3\sigma(\f{8n+20}{16})+\f83c(\f{8n+20}4)
\\&\q-2(\f43\sigma(\f{8n+20}{4})
-\f{16}3\sigma(\f{8n+20}{16})+\f83c(\f{8n+20}4))+\f43\sigma(2n+5)+\f83c(2n+5)
\\&=28\sigma(2n+5)-24\sigma(2n+5)-\f83\sigma(2n+5)-\f{16}{3}c(2n+5)+\f43\sigma(2n+5)
+\f83c(2n+5)
\\&=\f83(\sigma(2n+5)-c(2n+5)).\endalign$$
The proof is now complete. $\square$

 \pro{Theorem 4.3} Let $n\in\Bbb N$. Then
$$t(1,9,9,9;n)=\cases 16\sigma(\f{2n+7}9)&\t{if $n\e1\mod9,$}
\\0&\t{if $n\e2,4,5,7,8\mod9$},
\\\f 43(\sigma(2n+7)-c(2n+7))&\t{if $n\e0\mod 3$}.\endcases$$
\endpro
Proof. For $x\in\Bbb Z$ we see that $x(x-1)/2\e 0,1,3,6\mod 9$.
Thus, $t(1,9,9,9;n)=0$ for $n\e2,4,5,7,8\mod9$. Now we assume that
$n\e 0,1,3,6\mod 9$. For $n\e 1\mod 9$ we see that $9\mid 8n+28$ and
so
$$\align t(1,9,9,9;n)&=N_0(1,9,9,9;8n+28)
\\&=\big|\big\{(x,y,z,w)\in\Bbb Z^4\bigm|
8n+20=(3x)^2+9y^2+9z^2+9w^2,\ 2\nmid xyzw\big\}\big|
\\&=N_0\big(1,1,1,1;\f{8n+28}9\big)=t\big(1,1,1,1;\f{n-1}9\big)
=16\sigma\big(\f{2n+7}9\big).
 \endalign$$
For $n\e 0\mod 3$ we see that
$$\align t(1,1,9,9;n+1)
&=N_0(1,1,9,9;8n+28)
\\&=\big|\big\{(x,y,z,w)\in\Bbb Z^4\bigm|
8n+28=(3x)^2+y^2+9z^2+9w^2,\ 2\nmid xyzw\big\}\big|
\\&\q+\big|\big\{(x,y,z,w)\in\Bbb Z^4\bigm|
8n+28=x^2+(3y)^2+9z^2+9w^2,\ 2\nmid xyzw\big\}\big|
\\&=2N_0(1,9,9,9;8n+28)=2t(1,9,9,9;n).\endalign$$
Now applying Theorem 4.2 we deduce the result in this case.
$\square$

\pro{Theorem 4.4} Let $n\in\Bbb N$. Then
$$ t(1,1,1,9;n)=\cases 4\sigma(2n+3)+12\sigma(\f{2n+3}9)\
&\t{if $n\e0\mod 3$,}
\\8\sigma(2n+3)&\t{if $n\e1\mod 3$,}
\\4(\sigma(2n+3)-c(2n+3))&\t{if $n\e2\mod 3$.}\endcases$$
\endpro
Proof. For $n\e 0\mod 3$ we see that $3\mid 8n+12$. If
$8n+12=x^2+y^2+z^2+w^2$ for $x,y,z,w\in\Bbb Z$, then either $x\e y\e
z\e w\e 0\mod 3$ or $xyzw\e \pm 3\mod 9$. Thus,
$$N_0(1,1,1,1;8n+12)=4N_0(1,1,1,9;8n+12)-3N_0(9,9,9,9;8n+12).$$
This together with (4.1) yields
$$t(1,1,1,1;n+1)=\cases 4t(1,1,1,9;n)&\t{if $n\e 0,6\mod 9$,}
\\4t(1,1,1,9;n)-3t(1,1,1,1;\f{n-3}9)&\t{if $n\e 0\mod 3$.}
\endcases$$ Now combining the above with (1.2) yields the result in this
case.
\par Suppose $n\e 1\mod3$. Then $8n+12\e 2\mod 3$. If $8n+12
=x^2+y^2+z^2+9w^2$ for $x,y,z,w\in\Bbb Z$, then $3\mid xyz$ but
$9\nmid xyz$. Thus,
$$\align t(1,1,1,9;n)&=N_0(1,1,1,9;8n+12)
\\&=\big|\big\{(x,y,z,w)\in\Bbb Z^4\bigm|
8n+12=(3x)^2+y^2+z^2+9w^2,\ 2\nmid xyzw\big\}\big|
\\&\q+\big|\big\{(x,y,z,w)\in\Bbb Z^4\bigm|
8n+12=x^2+(3y)^2+z^2+9w^2,\ 2\nmid xyzw\big\}\big|
\\&\q+\big|\big\{(x,y,z,w)\in\Bbb Z^4\bigm| 8n+12=x^2+y^2+(3z)^2+9w^2,\
2\nmid xyzw\big\}\big|
\\&=3N_0(1,1,9,9;8n+12)
=3t(1,1,9,9;n-1).\endalign$$
This together with Theorem 4.2 yields
the result in this case.
\par For $n\e 2\mod 3$ we see that $8n+12\e 1\mod 3$ and so
$$\align t(1,1,1,9;n)&=N_0(1,1,1,9;8n+12)
\\&=\big|\big\{(x,y,z,w)\in\Bbb Z^4\bigm|
8n+12=(3x)^2+(3y)^2+z^2+9w^2,\ 2\nmid xyzw\big\}\big|
\\&\q+\big|\big\{(x,y,z,w)\in\Bbb Z^4\bigm|
8n+12=x^2+(3y)^2+(3z)^2+9w^2,\ 2\nmid xyzw\big\}\big|
\\&\q+\big|\big\{(x,y,z,w)\in\Bbb Z^4\bigm| 8n+12=(3x)^2+y^2+(3z)^2+9w^2,\
2\nmid xyzw\big\}\big|
\\&=3N_0(1,9,9,9;8n+12)=3t(1,9,9,9;n-2).\endalign$$
 Now combining the above with Theorem
4.3 yields the result in the case $n\e 2\mod 3$. The proof is now
complete. $\square$
 \par\q\par In conclusion we pose the following
conjecture.
 \pro{Conjecture 4.1} Suppose  $n\in\Bbb N$ and
$8n+9=3^{\beta}n_1$ with $3\nmid n_1$. Then
$$ t(1,1,3,4;n)=\f 12\Big(3^{\beta+1}\Ls 3{n_1}-1\Big)
\sum_{d\mid n_1}d\Ls 3d -\sum\Sb a,b\in\Bbb N,\ 2\nmid
a\\4(8n+9)=a^2+3b^2\endSb (-1)^{\f{a-1}2}a.$$
\endpro
\par Conjecture 4.1 has been checked for $n\le 1000$.

\par\q
\newline{\bf Acknowledgement}
\newline The second author is supported by the National Natural Science
Foundation of China (grant No. 11371163).

\end{document}